# RELIABILITY ANALYSIS OF SYSTEMS SUBJECT TO MUTUALLY DEPENDENT COMPETING FAILURE PROCESSES WITH CHANGING DEGRADATION RATE

NOOSHIN YOUSEFI, DAVID W. COIT

*Industrial and System Engineering Department, Rutgers Univeristy*
*Piscataway, NJ, 08854, USA*
*Email: no.yousefi@rutgers.edu, Coit@soe.rutgers.edu*

In this paper, a new reliability model has been developed for a single system degrading stochastically which experiences soft and hard failure. Soft failure occurs when the physical deterioration level of the system is greater than a predefined failure threshold, and hard failure occurs when the instantaneous stress caused by a shock process is greater than a critical threshold. It is considered that the degradation and shock process are mutually dependent. In fact, each arriving shock accelerates the degradation process by adding abrupt additional damages to the degradation path and changing the degradation rate according to specific magnitude; also, the cumulative degradation changes the occurrence intensity of shock process. A gamma process is used as a stochastic process to model the degradation path. A realistic numerical example is presented to illustrate the proposed reliability.

## 1. Introduction and Background

Industrial systems and products experience multi-failure modes. A degradation process and shock process are the two most common failure modes in many systems. Degradation processes present the level of system physical deterioration over time. When system degradation reaches a pre-determined critical value, the system fails due to soft failure. Another competing failure process is called a hard failure, which is due to instantaneous stress caused by shock process. There has been significant research considering soft failure and hard failure in system reliability [1-2].

In recent research, two types of dependencies between shock and degradation processes have been studied; shock-degradation dependence and degradation-shock dependence models. In shock-degradation dependence models, random shock has damages on the degradation process, such as causing an abrupt increase on the degradation process and increasing the degradation rate. Peng et al. [3] developed a reliability model for a system when arriving shocks cause an additional abrupt change to the cumulative degradation process. Song et al [4] studied the reliability and condition-based maintenance of a multi-component system considering shock-degradation dependence. Rafiee et al [5] investigated a reliability model considering changing of degradation rate by different shock patterns.

The majority of recent research focused on shock-degradation models, while in some systems shocks can be influenced by the degradation process. There has been significant research on degradation-shock models such as Fan et al [6], Zhu et al [7] and Ye et al [8]

However, in practice, systems can experience mutual dependency of shock and degradation processes. While a shock process can be influenced by the degradation process and it can also have effects on degradation process. Che et al. [9] proposed a reliability model considering mutually dependent shock and degradation processes. Each shock has an additional abrupt change on the degradation process and the cumulative degradation changes the occurrence intensity of coming shocks. They modeled the shock process using facilitation model and the degradation process follows a deterministic degradation model.

In this paper, mutually dependence of shock and degradation process is considered for a system where coming shocks not only have additional abrupt on degradation path, but also change the degradation rate. Some systems deteriorate faster when they become more prone to failure; so, beside the abrupt damage of each shock exposure on the degradation process, they can change the degradation rate to make it deteriorate faster. In addition, it is assumed that current cumulative degradation can increase the occurrence intensity of shock process. In this paper, to consider the uncertainty of degradation over time, a gamma process is used as a stochastic process, which can model the degradation path of systems degrading in a form of cumulative damage.

## 2. System reliability analysis

Two dependent competing failure processes are considered: soft and hard failure. It is considered that arrival shocks have two types of damage on the cumulative degradation process. First, each incoming shock has an additional abrupt change on the cumulative degradation process. Second, if the shock magnitude is greater than a predefined critical value, the degradation rate changes from the arrival time of that shock. Moreover, the cumulative degradation increases the occurrence intensity of the next shocks. The cumulative degradation is the summation of pure degradation and additional abrupt changes caused by arriving shocks. Therefore, it can be concluded that next shock arrival is influenced by degradation and the total number of arrived shocks, and it has damage effects on the cumulative degradation.

A Poisson process is a random process when each event is independent of all other events in the process. However, in this study the arrival shock depends on the number of shocks arriving to the system and the total degradation, so the Poisson process is not suitable for the shock arrival process. Che et al [9] showed that random shock process $\{N(t), t>0\}$ can be modeled as a facilitation model using intensity function $\lambda_i(t)$. As it is proved in [9] the probability of having $i$ shocks by time $t$ can be calculated as the following equation

$$P_i(t) = \frac{(\eta^{-1}+(i-1))^i (\exp(\eta\Lambda_0(t))-1)^i (\exp(-(1+\eta i)\Lambda_0(t)))}{i!} \qquad (1)$$
$$= C_{\eta^{-1}+i-1}^{i}(1-\exp(-\eta\Lambda_0(t)))^i (\exp(-\eta\Lambda_0(t)))^{\eta^{-1}}$$

Where $\Lambda_i(t) = \int_0^t \lambda_i(t)$ and the intensity function is $\lambda_i(t) = (1+\eta i)\lambda_0(t)$ which depends on the total degradation $\lambda_0(t) = \lambda_0 + \gamma X_S(t)$. $\gamma$ indicates the effect of the current degradation on the intensity and $\eta$ is the facilitation factor which shows the effect of abrupt degradation on the shock process, and total degradation is the summation of degradation process and damage caused by shocks $X_S(t) = X(t,\theta) + \sum_{i=1}^{N(t)} Y_i$.

To calculate the system reliability by time $t$, the system should not experience any hard and soft failure by time $t$. It is assumed that the shock damage is an $i.i.d$ random variable which follows a normal distribution $Y_{ij} \sim Normal(\mu_{Yi}, \sigma^2_{Yi})$. The probability of having no soft failure by time $t$ is ($P_{NS}$)

$$P_{NS}(t) = P(X_S(t) < H) = \sum_{m=0}^{\infty} P(X(t,\theta) + \sum_{j=1}^{N(t)} Y_j < H) \times P(N(t) = m) \qquad (2)$$

It is assumed that the degradation can be modeled as gamma process. For $t>0$ and $s>0$, $X(t) - X(s) \sim gamma(\alpha(t) - \alpha(s), \beta)$. Where $\alpha(t)$ is the shape parameter and $\beta$ is the scale parameter.

$$g(\theta; \alpha(t-s), \beta) = \frac{\beta^{\alpha(t-s)} \theta^{\alpha(t-s)-1} \exp(-\beta\theta)}{\Gamma(\alpha(t-s))} \qquad (3)$$

Moreover, it is assumed that the shock magnitude is an $i.i.d$ random variable follows normal distribution $W_{ij} \sim Normal(\mu_{Wi}, \sigma^2_{Wi})$, so the probability of no hard failure by time $t$ is as follow ($P_{NH}$)

$$P_{NH}(t) = P(W_j < D) = F_W(D) = \Phi(\frac{D - \mu_W}{\sigma_W}) \quad (4)$$

The system reliability by time *t* can be calculated using equation (5)

$$R(t) = \sum_{m=0}^{\infty} \left[ (F_W(D))^m P\left( X(t) + \sum_{j=1}^{N(t)} Y_j < H \right) | N(t) = m \right] P(N(t) = m)$$

$$= \sum_{m=0}^{\infty} \left[ (F_W(D))^m \times \int_0^H G(H - y; \alpha t, \beta) f_Y^{<m>}(y) dy \right] P(N(t) = m) \quad (5)$$

Using equation (1) the probability that we have *m* shocks by time *t* is as follow.

$$P(N(t) = m) = \int_\theta \int_Y P(N(t) = m | \sum_{i=1}^{N(t)} Y_i = y, \theta = z) \times f_\theta(z) dz f_Y^{<m>}(y) dy$$

$$= \int_0^{H-y} \int_0^H C_{\eta^{-1}+m-1}^m (1 - \exp(-\eta \Lambda_0(t)))^m (\exp(-\eta \Lambda_0(t)))^{\eta^{-1}} \times f_\theta(z) dz f_Y^{<m>}(y) dy$$

$$= \int_0^{H-y} \int_0^H C_{\eta^{-1}+m-1}^m \left( 1 - \exp\left( -\eta \times \int_0^t \lambda_0 + \gamma(X(v; z) + y)) dv \right) \right)^m \left( \exp\left( -\eta \int_0^t \lambda_0 + \gamma(X(v; z) + y)) dv \right) \right)^{\eta^{-1}} \times f_\theta(z) dz f_Y^{<m>}(y) dy$$

(6)

Where $f_\theta(z)$ is the probability density function of $\theta$. By substituting equation (6) in equation (5) the system reliability is as follow:

$$R(t) = \sum_{m=0}^{\infty} \left[ P(W_j < D)^m P\left( X(t) + \sum_{j=1}^{N(t)} Y_j < H \right) | N(t) = m \right] P(N(t) = m)$$

$$= \sum_{m=0}^{\infty} \left[ P(W_j < D)^m \times \int_0^H G(H - y; \alpha t, \beta) f_Y^{<m>}(y) dy \right]$$

$$\times \int_0^{H-y} \int_0^H C_{\eta^{-1}+m-1}^m \left( 1 - \exp\left( -\eta \times \int_0^t \lambda_0 + \gamma(X(v; z) + y)) dv \right) \right)^m \left( \exp\left( -\eta \int_0^t \lambda_0 + \gamma(X(v; z) + y)) dv \right) \right)^{\eta^{-1}} \times f_\theta(z) dz f_Y^{<m>}(y) dy$$

(7)

## 3. Reliability analysis considering changing degradation rate

In this paper, it is considered that the degradation rate changes if any shock magnitude is greater than a predefined damage threshold. Therefore, two thresholds should be defined for shock magnitude, one is for detecting hard failure ($D_1$), and the other for damage which is in form of changing degradation rate ($D_0$). Figure 1 shows this model. At time $t_3$ the shock magnitude is greater than the hard failure threshold which caused the system fails, and at time $t_2$ the shock magnitude is greater than the damage threshold which causes change of the degradation rate.

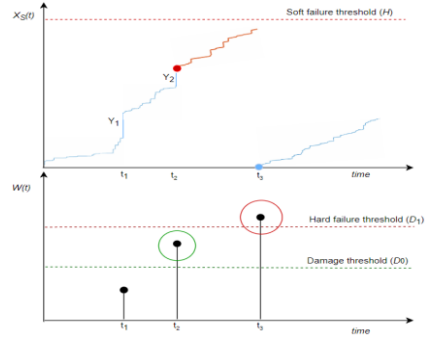

Figure 1. Shock model

To calculate the system reliability considering changing degradation rate, the probability of all the possible scenarios should be calculated and added together. There may have three different scenarios for this system.

1. When there is no shock has arrived by time $t$.

$$R_1(t) = P(X < H) \times P(N(t) = 0) = G(H; \alpha t, \beta) \times P(N(t) = 0)$$
$$= G(H; \alpha t, \beta) \times \int_0^H \exp\left(-\eta \int_0^t (\lambda_0 + \gamma X(v; z)) dv\right)^{\eta-1} f_\theta(z) dz \qquad (8)$$

2. When $m$ shocks have arrived but all of them are less than the damage threshold ($D_0$)

$$R_2(t) = \sum_{m=0}^{\infty} \left[ (F_W(D_0))^m \times \int_0^H G(H-y; \alpha t, \beta) f_Y^{<m>}(y) dy \right] \times P(N(t) = m)$$
$$= \sum_{m=0}^{\infty} \left[ (F_W(D_0))^m \times \int_0^H G(H-y; \alpha t, \beta) f_Y^{<m>}(y) dy \right]$$
$$\times \int_0^{H-y} \int_0^H C_{\eta^{-1}+m-1}^m \left(1 - \exp\left(-\eta \times \int_0^t \lambda_0 + \gamma(X(v; z) + y)) dv\right)\right)^m \left(\exp\left(-\eta \int_0^t \lambda_0 + \gamma(X(v; z) + y)) dv\right)\right)^{\eta-1} \times f_\theta(z) dz f_Y^{<m>}(y) dy$$
(9)

3. When $m$ shocks have arrived and $j^{\text{th}}$ random shock is greater than the damage threshold ($D_0$)

$$R_3(t) = \sum_{m=1}^{\infty} \sum_{j=1}^{m} \left[ (F_W(D_0))^{j-1} (F_W(D_1) - F_W(D_0))(F_W(D_1))^{m-j} \right]$$
$$\times \left[ \int_0^t P\left( X_1(t_j) + X_2(t - t_j) + \sum_{j=1}^{N(t)} Y_j < H \right) f_{t_j}(t) dt_j \right] \times P(N(t) - N(t_j) = m - j) P(N(t_j) = j)$$
$$= \sum_{m=1}^{\infty} \sum_{j=1}^{m} \left[ (F_W(D_0))^{j-1} (F_W(D_1) - F_W(D_0))(F_W(D_1))^{m-j} \right]$$
$$\times \int_0^H \int_0^t G(H - y; \alpha_1(t_j) + \alpha_2(t - t_j), \beta) \times P(N(t) - N(t_j) = m - j) P(N(t_j) = j) f_{t_j}(t) dt_j f_Y^{<m>}(y) dy$$

Where the probability density function of $t_j$ can be calculated as $f_{t_j}(t) = \dfrac{dF_{t_j}(t)}{dt}$ and

$F_{t_j}(t) = P(T \leq t_j) = 1 - P(T > t_j) = 1 - P_{t_j}(X \leq j-1)$. The probability of having $j$ shocks by time $t_j$ should be calculated using equation (6) while the shape parameter for gamma process is $\alpha_1(t)$ and scale parameter is $\beta$. To calculate the probability of having $m$-$j$ shocks from $t_j$ to $t$, equation (6) can be used, while the initial degradation is the cumulative degradation up to time $t_j$.

$\Lambda_0(t-t_j) = \int_0^H \int_0^{t-t_j} \lambda_0 + \gamma\left(X(v) + y + \left(\sum_{k=0}^{j} Y_k + X(t_j)\right)\right) f_Y^{<m-j>}(y) dy$. $X(t)$ follows a gamma process with shape parameter $\alpha_2(t)$ and scale parameter $\beta$.

## 4. Numerical results

A jet pipe servo valve is considered as a realistic numerical example in this section. The valve control hydraulic oil flows by the spool sliding in the sleeve [10]. However, the oil pollution can lead to the breakdown of leakage and the contamination lock (i.e. the clamping stagnation failure [11]). In this paper, the contamination lock is considered as hard failure and the wear between spool and sleeve is considered as soft failure. Table 1 shows the parameters for reliability analysis of this example.

Table 1 Parameters value

| Parameters | Description | Valve | Source |
|---|---|---|---|
| $H$ | The soft failure threshold | 5 mm | Fan et al. [12] |
| $D_1$ | The Hard failure threshold | 40 N | Che et al. [9] |
| $D_0$ | The damage threshold | 30 N | Assumption |
| $\alpha_1$ | Initial Shape parameter of gamma process | 0.5 | Assumption |
| $\alpha_2$ | Changed Shape parameter of gamma process | 0.9 | Assumption |
| $\beta$ | Scale parameter of gamma process | 1.2 | Assumption |
| $\lambda_0$ | The initial intensity of random shock | $2.5 \times 10^{-5}$ | Fan et al. [12] |
| $\eta$ | The facilitation factor | 0.2 | Che et al. [9] |
| $\gamma$ | The dependence factor | 0.001 | Fan et al. [12] |
| $W$ | The friction caused by contamination lock | $W \sim N(10, 5^2)$ N | Assumption |
| $Y$ | The wear increase by contamination lock | $Y \sim N(0.5, 0.1^2)$ mm | Che et al. [9] |

Monte Carlo simulation with $10^5$ replications is used to calculate the system reliability. Figure 2 shows the system reliability with and without considering changing degradation rate. When the degradation rate changes, the system degrades faster, and subsequently there are more shocks coming to the system, so as it is shown in Figure 2, the reliability considering changing degradation rate is lower than reliability with fixed degradation rate.

A sensitivity analysis is performed for effect of the damage threshold on the system reliability. As it is shown in Figure 3, as damage threshold decreases, the system reliability decreases faster. When the damage threshold is high, there is fewer changes in the system degradation rate. While it is low, we will have more shocks greater than the damage threshold, and consequently the degradation rate changes to make the system degrade faster.

Figure 4 shows the sensitivity analysis of system reliability on degradation-shock dependence factor ($\gamma$). When $\gamma$ is small, there is less degradation-shock dependence, and high $\gamma$ makes the system receive more shocks which results in more shocks with magnitude of greater than the damage threshold, and subsequently faster deterioration. As it is shown on Figure 4, the system with high $\gamma$ degrades faster and has lower system reliability, while there is almost no degradation-shock dependence when $\gamma$ is high, so the system reliability is higher.

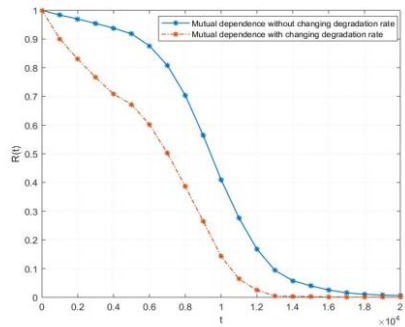

Figure 2 System reliability considering changing degradation rate

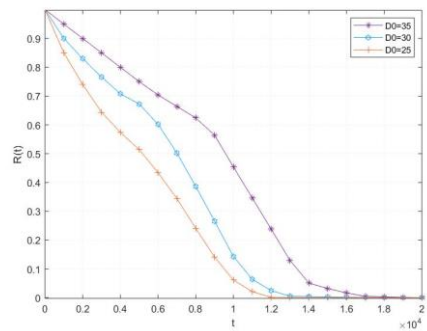

Figure 3 Sensitivity analysis on damage threshold

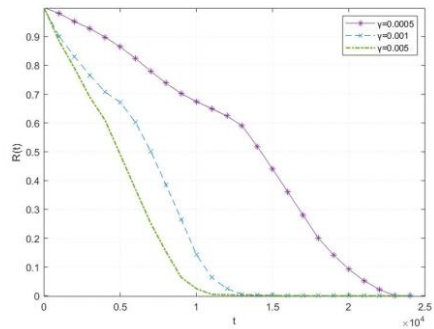

Figure 4 Sensitivity analysis of $R(t)$ on $\gamma$

## 5. Conclusion

In this paper, a new reliability model is developed for a single component system degrading as a gamma process. The system experiences two competing failure processes of soft and hard failure. If the system's total degradation, which is summation of pure degradation and additional abrupt damage caused by shocks, is greater than a predefined threshold, the soft failure happens, and when any shock magnitude is greater than a critical threshold, the system experiences hard failure. In this paper, the mutual dependency of shock and degradation processes are considered. Each incoming shock has two types of damage on the system degradation. It has additional jumps on the degradation path and it can change the degradation rate, if its magnitude is greater than a damage threshold. Consequently, the cumulative degradation changes the occurrence intensity of next shocks. A realistic example of a jet pipe servo valve is used to show the performance of the new reliability model compared to previous models.